\documentclass[twoside,leqno,11pt]{article}

\usepackage{graphics,ifthen}
\usepackage{latexsym}
\usepackage{amsmath}
\usepackage{amssymb}
\usepackage{mathrsfs}
\begin{document}
\input{latexP.sty}
\input{referencesP.sty}
\input epsf.sty

\def\ind{\stackrel{\mathrm{ind}}{\sim}}
\def\iid{\stackrel{\mathrm{iid}}{\sim}}
\def\Prodi{\mathop{{\lower9pt\hbox{\epsfxsize=15pt\epsfbox{pi.ps}}}}}
\def\prodi{\mathop{{\lower3pt\hbox{\epsfxsize=7pt\epsfbox{pi.ps}}}}}

\def\Definition{\stepcounter{definitionN}\
    \Demo{Definition\hskip\smallindent\thedefinitionN}}
\def\EndDefinition{\EndDemo}
\def\Example#1{\Demo{Example [{\rm #1}]}}
\def\EndExample{\qed\EndDemo}
\def\Category#1{\centerline{\Heading #1}\rm}
\
%% Paper specific definitions
\def\e{\text{\hskip1.5pt e}}
\newcommand{\eps}{\epsilon}
\newcommand{\proof}{\noindent {\bf Proof:\ }}
\newcommand{\remarks}{\noindent {\bf Remarks:\ }}
\newcommand{\note}{\noindent {\bf Note:\ }}
\newcommand{\examp}{\noindent {\bf Example:\ }}
\newcommand{\Lower}[2]{\smash{\lower #1 \hbox{#2}}}
\newcommand{\ben}{\begin{enumerate}}
\newcommand{\een}{\end{enumerate}}
\newcommand{\bi}{\begin{itemize}}
\newcommand{\ei}{\end{itemize}}
\newcommand{\hp}{\hspace{.2in}}

\newtheorem{lw}{Proposition 3.1, Lo and Weng (1989)}
\newtheorem{thm}{Theorem}[section]
\newtheorem{defin}{Definition}[section]
\newtheorem{prop}{Proposition}[section]
\newtheorem{lem}{Lemma}[section]
\newtheorem{cor}{Corollary}[section]
\newcommand{\rb}[1]{\raisebox{1.5ex}[0pt]{#1}}
\newcommand{\mc}{\multicolumn}
%Mathrsfs Font
\newcommand{\Bcr}{\mathscr{B}}
\newcommand{\Ucr}{\mathscr{U}}
\newcommand{\Gcr}{\mathscr{G}}
\newcommand{\Dcr}{\mathscr{D}}
\newcommand{\CS}{\mathscr{C}}
\newcommand{\Fcr}{\mathscr{F}}
\newcommand{\Icr}{\mathscr{I}}
\newcommand{\Lcr}{\mathscr{L}}
\newcommand{\Mcr}{\mathscr{M}}
\newcommand{\Ncr}{\mathscr{N}}
\newcommand{\Pcr}{\mathscr{P}}
\newcommand{\Qcr}{\mathscr{Q}}
\newcommand{\Scr}{\mathscr{S}}
\newcommand{\Tcr}{\mathscr{T}}
\newcommand{\Vcr}{\mathscr{V}}
\newcommand{\Xcr}{\mathscr{X}}
\newcommand{\Ycr}{\mathscr{Y}}
%Mathbb Font
\newcommand{\E}{\mathbb{E}}
\newcommand{\F}{\mathbb{F}}
\newcommand{\I}{\mathbb{I}}
\newcommand{\Q}{\mathbb{Q}}
\newcommand{\X}{\mathbb{X}}
\newcommand{\Pe}{\mathbb{P}}
\newcommand{\M}{\mathbb{M}}
\newcommand{\Wbb}{\mathbb{W}}

\def\Beta{\text{Beta}}
\def\Dir{\text{Dirichlet}}
\def\DP{\text{DP}}
\def\P{{\bf p}}
\def\fhat{\widehat{f}}
\def\GA{\text{gamma}}
\def\ind{\stackrel{\mathrm{ind}}{\sim}}
\def\iid{\stackrel{\mathrm{iid}}{\sim}}
\def\J{{\bf J}}
\def\K{{\bf K}}
\def\min{\text{min}}
\def\N{\text{N}}
\def\m{{\bf m}}
\def\p{{\bf p}}
\def\U{{\bf U}}
\def\v{{\bf v}}
\def\V{{\bf V}}
\def\W{{\bf W}}
\def\w{{\bf w}}
\def\T{{\bf T}}
\def\m{{\bf m}}
\def\MM{{\bf M}}
\def\X{{\bf X}}
\def\y{{\bf y}}
\def\Y{{\bf Y}}
\def\tps{\mbox{\scriptsize ${\theta H}$}}   %   smaller "\psi"-vector
\def\ups{\mbox{\scriptsize ${P_{\theta}(g)}$}}   %   smaller "\psi"-vector
\def\vps{\mbox{\scriptsize ${\theta}$}}   %   smaller "\psi"-vector
\def\vups{\mbox{\scriptsize ${\theta >0}$}}   %   smaller "\psi"-vector
\def\hps{\mbox{\scriptsize ${H}$}}   %   smaller "\psi"-vector
\def\rps{\mbox{\scriptsize ${(\theta+1/2,\theta+1/2)}$}}   %   smaller "\psi"-vector
\def\sps{\mbox{\scriptsize ${(1/2,1/2)}$}}   %   smaller "\psi"-vector

\newcommand{\reals}{{\rm I\!R}}
\newcommand{\PR}{{\rm I\!P}}
\def\Z{{\bf Z}}
\def\yy{{\mathcal Y}}
\def\rr{{\mathcal R}}
\def\BP{\text{beta}}
\def\ts{\tilde{t}}
\def\js{\tilde{J}}
\def\gs{\tilde{g}}
\def\fs{\tilde{f}}
\def\ys{\tilde{Y}}
\def\ps{\tilde{\mathcal {P}}}

\def\Report{Lancelot F. James and John W. Lau}
\def\Author{Hyperbolic Integrated Volatility}
\pagestyle{myheadings}
\markboth{\Author}{\Report}
\thispagestyle{empty}

\bct\Heading A Class of Generalized Hyperbolic Continuous Time\\
Integrated Stochastic Volatility  Likelihood Models\lbk\lbk\smc
Lancelot F. James and John W. Lau\footnote{ \eightit AMS 2000
subject classifications.
               \rm Primary 62G05; secondary 62F15.\\
\eightit Corresponding authors address.
                \rm The Hong Kong University of Science and Technology,
Department of Information and Systems Management, Clear Water Bay,
Kowloon, Hong Kong.
\rm lancelot\at ust.hk\\
\indent\eightit Keywords and phrases.
                \rm
          Black-Scholes models,
          Chinese Restaurant Process,
          Integrated stochastic volatility,
          Latent GARCH,
          Ornstein-Uhlenbeck process,
          Poisson random measure.
          }
\lbk\lbk \BigSlant The Hong Kong University of Science and
Technology and University of Bristol\rm \lbk %(\today)%
\ect \Quote This paper discusses and analyzes a class of
likelihood models which are based on two distributional
innovations in financial models for stock returns. That is, the
notion that the marginal distribution of aggregate returns of
log-stock prices are well approximated by generalized hyperbolic
distributions, and that volatility clustering can be handled by
specifying the integrated volatility as a random process such as
that proposed in a recent series of papers by Barndorff-Nielsen
and Shephard~(BNS). Indeed, the use of just the integrated
Ornstein-Uhlenbeck(INT-OU) models of BNS serves to handle both
features mentioned above. The BNS models produce likelihoods for
aggregate returns which can be viewed as a subclass of latent
regression models where one has n conditionally independent Normal
random variables whose mean and variance are representable as
linear functionals of a common unobserved Poisson random measure.
James~(2005b) recently obtains an exact analysis for such models
yielding expressions of the likelihood in terms of quite tractable
Fourier-Cosine integrals. Here, our idea is to analyze a class of
likelihoods, which can be used for similar purposes, but where the
latent regression models are based on n conditionally independent
models with distributions belonging to a subclass of the
generalized hyperbolic distributions and whose corresponding
parameters are representable as linear functionals of a common
unobserved Poisson random measure. Our models are perhaps most
closely related to the Normal inverse Gaussian/GARCH/A-PARCH
models of Brandorff-Nielsen~(1997) and Jensen and Lunde~(2001),
where in our case the GARCH component is replaced by quantities
such as INT-OU processes. It is seen that, importantly, such
likelihood models exhibit quite different features structurally.
Rather than Fourier-Cosine integrals, the exact analysis of these
models yields characterizations in terms of random partitions of
the integers which can be easily handled by Bayesian SIS/MCMC
procedures similar to those which have been applied to
Dirichlet/Gamma process mixture models. Importantly, these methods
do not necessarily require the simulation of random measures. One
nice feature of the model is that it allows for more flexibility
in terms of modelling of external regression parameters. Our
models may also be viewed as alternatives to closely related
latent class GARCH models arising in financial economics.
 \EndQuote
%\baselineskip14pt
%\begin{document}
\rm
%\newpage

\section{Introduction}
In financial economics, it is well known that Gaussian based
models such as the Black-Scholes-Samuelson model for the log-stock
prices of returns do not fit with empirical observations when
returns are observed over moderately sized intervals. For example,
on a daily basis. The Black-Scholes-Samuelson model may be
described in terms of the following stochastic differential
equation, \Eq dx^{*}(t)=(\mu+\beta \sigma^{2})dt+\sigma dw(t)
\label{Black}\EndEq where $x^*(t)$ denotes the price level,
$\sigma^{2}$ represents a constant volatility and $w(t)$ is
Brownian Motion. When observed over $i=1,\ldots, n$ equally spaced
time intervals of length $\Delta>0$, one has that the aggregate
returns $x^{*}(i\Delta)-x^{*}((i-1)\Delta)$ are iid Normal random
variables with mean and variance
$(\mu\Delta+\sigma^{2}\beta,\sigma^{2})$. In terms of statistical
inference, this produces a classical Normal likelihood where
estimation of parameters $(\beta,\mu)$ are straightforward.
However, it is known that, while the model ~\mref{Black} is
plausible for large $\Delta$, when $\Delta$ is of moderate size
the aggregate returns exhibit behavior more like that of
semi-heavy tailed distributions. Moreover, these models exhibit a
feature known as {\it volatility persistence} or {\it clustering}.
This suggests that $\sigma^{2}$ should be replaced by a dynamic
random process which has correlated increments. See for instance
Carr and Wu~(2004), Carr,Geman, Madan and Yor~(2003),
Barndorff-Nielsen and Shephard~(2001a,b), Duan~(1995), and
Engle~(1982) for these points and various proposals to
enhance~\mref{Black}. Here we shall focus on the model of
Barndorff-Nielsen and Shephard~(2001a, b) which we now describe.
\subsection{BNS model and likelihood}
A quite attractive model was introduced by Barndorff-Nielsen and
Shephard~(2001a, b). Their proposed continuous time stochastic
volatility (SV) model is based on the following differential
equation, \Eq dx^{*}(t)=(\mu+\beta v(t))dt+v^{1/2}(t)dw(t)
\label{BNS}\EndEq where $x^*(t)$ denotes the price level, and
$v(t)$ is a stationary Ornstein-Uhlenbeck (OU) process which
models the {\it instantaneous volatility} and is independent of
$w(t)$. The induced likelihood model, which is based on the
integrated volatility $\tau(t)=\int_{0}^{t}v(u)du$, can be
described as follows. Let
${X_{i}}:=x^{*}(i\Delta)-x^{*}((i-1)\Delta),$ for $i=1,\ldots, n$
denote a sequence of the returns of the log price of a stock
observed over intervals of length $\Delta>0$. Additionally for
each interval $[(i-1)\Delta,i\Delta]$, let
$\tau_{i}=\tau(i\Delta)-\tau((i-1)\Delta)$. Now the model in
~\mref{BNS} implies that $X_{i}|\tau_{i},\beta, \mu$ are
conditionally independent with \Eq
X_{i}=\mu\Delta+\tau_{i}\beta+\tau^{1/2}_{i}\epsilon_{i}.
\label{data}\EndEq where $\epsilon_{i}$ are independent standard
Normal random variables. Hence if $\tau$ depends on external
parameters $\theta$, one is interested in estimating
$(\mu,\beta,\theta)$ based on the likelihood \Eq
\Lcr_{BNS}(\X|\mu,\beta,\theta)=\int_{{\mathbb R}^{n}_{+}}
\[\prod_{i=1}^{n}\phi(X_{i}|\mu\Delta+\beta
\tau_{i},\tau_{i})\]f(\tau_{1},\ldots,\tau_{n}|\theta)d\tau_{1},\ldots,d\tau_{n}
\label{BNSlik} \EndEq where, setting $A_{i}=(X_{i}-\mu\Delta)$,
and ${\bar A}=n^{-1}\sum_{i=1}^{n}A_{i}$,
$$\phi(X_{i}|\mu\Delta+\beta
\tau_{i},\tau_{i})={\mbox e}^{A_{i}\beta
}\frac{1}{\sqrt{2\pi}}\tau^{-1/2}_{i}{\mbox
e}^{-A^{2}_{i}/(2\tau_{i})}{\mbox e}^{-\tau_{i}\beta^{2}/2}
$$
denotes a Normal density. We note that because of the complex
dependence structure of the joint density
$f(\tau_{1},\ldots,\tau_{n}|\theta)$, the likelihood was thought
to be intractable. Thus inhibiting full likelihood based
statistical inference, for models involving quite arbitrary
$\tau$.  This is in contrast to the case of {\it latent class}
GARCH models which are of considerable interest in financial
economics[see for instance, Fiorentini, Sentana, and
Shephard~(2004)]. However, in a closely related recent paper
James~(2005b) shows that the likelihood~\mref{BNSlik}, where the
$\tau_{i}$ are further generalized to be linear functionals of a
Poisson random measure, is tractable and can be expressed exactly
in terms of multi-dimensional Fourier-cosine transforms. The
implication is that in general one can use classical numerical
techniques to evaluate the likelihood. Moreover, these expressions
are similar to quantities which regularly appear in the math
finance literature on option pricing and related areas.

In this paper we offer another approach that still allows us to
work with integrated OU processes and indeed more general objects.
Our purpose is two fold. One to propose models which we believe
are complementary to the above framework but exhibit quite
different features, and in fact more flexibility in the sense of
incorporating more general regression coefficients. Secondly, we
believe that because these two models exhibit different features
that this invites individuals of quite varying backgrounds to
conduct research on similar topics. Of course, exact analysis of
these two classes of models also allows one to more easily
critique, compare  and improve such models. One could also
consider a third class of models based on a hybridization of the
two models.

The difference is that our approach yields expressions of the
likelihood in terms of random partitions of the integers which can
considered relatives of the Blackwell and MacQueen~(1973) P\'olya
urn distribution. Hence these models inherit many of the
well-known features of Dirichlet/Gamma process mixture models and
extensions addressed in James~(2005a). Additionally, the posterior
distribution of the random processes are also more in line with
what happens for the case of Bayesian multiplicative intensity
models, as it depends on the jumps of the underlying Poisson
random measure. Moreover our models serve as an alternative to
{\it latent class} GARCH models. In fact, one will see
 in the next section that our models are perhaps most closely related to the Normal inverse
Gaussian/GARCH/A-PARCH models of Brandorff-Nielsen~(1997) and
Jensen and Lunde~(2001), where we replace their GARCH/A-PARCH
components by general $\tau$.

\section{A class of generalised hyperbolic integrated stochastic
volatility models} Before we present the model, we note that many
authors have fitted semi-heavy tailed models in finance by
specifying $\sigma^{2}$ in the Black-Scholes formula to be a {\it
generalized inverse Gaussian} (GIG) distribution.  Hence the
aggregate returns are from a {\it generalised Hyperbolic}~(GH)
distribution. We pause to describe this density which we shall use
later. Let $\lambda$, $v$ and $\delta$ be such that
$-\infty<\lambda<\infty$, while $v$ and $\delta$ are non-negative
and not simultaneously $0$. As in Barndorff-Nielsen and
  Shephard(2001a), $T$ is GIG$(\lambda,
  \delta, v)$ random variable if its density is of the form
  $$
  f_{GIG}(t|\lambda,\delta,v)=\frac{{({v/\delta})}^{\lambda}}{2K_{\lambda}(\delta
  v)}t^{\lambda-1}\exp\{-\frac{1}{2}(\delta^{2}t^{-1}+v^{2}t)\}
  $$
  where $K_{\lambda}$ is a Bessel function. When $\delta=0$ and $\lambda>0$, $v>0$ , GIG$(\lambda,0,v)$
equates with the Gamma distribution.  When $\lambda<0$, $\delta>0$
and $v=0$, then GIG$(\lambda, \delta,0)$ is a reciprocal, or
inverse Gamma distribution. Using the parametrization,
$\lambda=-a$, for $a>0$, and $b=\delta^{2}/2$, yields the density
of an inverse Gamma distribution with parameters, $a,b$. A special
case of this is when $\lambda=-1/2$ leading to a stable law of
index $1/2$. The inverse Gaussian distribution defined by setting
$\lambda=-1/2$,$\delta>0$, and $v>0$ that is a
GIG$(-1/2,\delta,v)$. The Hyperbolic distribution coincides with
the case of $\lambda=1$.  See Prause~(1999) and Eberlein~(2001)
for some additional background and references. The additional
innovation in, for instance,  Barndorff-Nielsen and
  Shephard(2001a, b) is that modelling volatility as a random process, $v(t)$, rather than
a random variable, not only allows for semi-heavy-tailed models,
but additionally induces serial dependence.
\subsection{The model and conditional likelihood}
We now describe a model which is a direct variant of~\mref{data}
but is otherwise a subclass of a considerably more flexible but
still tractable proposal which is given in section 2.3. Let
$X_{s,t}:=x^{*}(t)-x^{*}(s)$ denote the aggregate return of the
log stock price over some interval$[s,t]$ for $0\leq s<t$.
Furthermore, define $\tau_{s,t}=\tau(t)-\tau(s)$ and
$A_{s,t}=(x-\mu(t-s))$. Then given some filtration
$\Fcr_{\tau(t)}$ determined by $\tau$, and further depending on
$\beta$ and $\mu$, we assume that $X_{s,t}$ is conditionally
independent of the past with density \Eq
f_{X_{s,t}}(x|\tau,\beta,\mu)=
{\(\frac{2}{\pi}\)}^{1/2}\frac{\tau^{\lambda}_{s,t}{\mbox
e}^{A_{s,t}\beta}}{\Gamma(\lambda)}{\(\sqrt{\frac
{\beta^{2}}{[2\tau_{s,t}+{A_{s,t}}^{2}]}}\)}^{\lambda+1/2}K_{\lambda+1/2}
(\beta\sqrt{(2\tau_{s,t}+{A_{s,t}}^{2})}) \label{GHSV}\EndEq for
$\lambda>0.$

It follows that the density $~\mref{GHSV}$, for fixed $\tau$,
represents a subclass of generalized Hyperbolic~(GH) densities
which reduces to the Student distribution when $\mu=0$ and
$\beta=0$, but otherwise is one of the well-defined limiting cases
of the (GH) model[see for instance Prause~(1999)]. We point out
further that although this density, for fixed $\tau$, does not
contain the Normal inverse Gaussian or Hyperbolic distribution,
Prause~(1999, p. 7-11) gives examples where the densities
in~\mref{GHSV} provides a more plausible fit to the data than
those models. However, due to the general distributional
flexibility of $\tau$, these issues do not really concern us, and
we shall further take $\lambda=1$ for additional tractability.
That is to say estimation and model fitting will depend on
parameters such as $\theta,\beta$ and $\mu$ and the distributional
features of $\tau$.

It would appear that the models in~\mref{GHSV}, and its
corresponding likelihood model, are considerably more complex
than~\mref{BNS},~\mref{data} and~\mref{BNSlik}. However for each
increment one may write, \Eq X_{s,t}\overset
{d}=\mu(t-s)+(\tau_{s,t}/Z)\beta+{(\tau_{s,t}/Z)}^{1/2}\epsilon_{s,t}
\label{genrv}\EndEq where $Z$ is a Gamma random variable with
shape $\lambda$, which we set to $1$, and scale $1$, and
$\epsilon_{s,t}$ is an independent standard Normal random
variable. Moreover, by utilizing a change of variable
$W=Z/\tau_{s,t}$, the density~\mref{GHSV} can be written as
\Eq\int_{0}^{\infty}\phi(x|\mu(t-s)+\beta w^{-1},w^{-1}){\mbox
e}^{-w\tau_{s,t}}\tau_{s,t}dw\label{mixrep}\EndEq
\subsection{The likelihood}
Let $Z_{i}$ be iid Gamma random variables with shape $\lambda=1$
and scale $1$. Then under the setting of~\mref{BNSlik} our model
translates into the case where $X_{i}|\tau_{i},\mu,\beta$ are
conditionally independent and representable as \Eq X_{i}\overset
{d}=\mu\Delta+(\tau_{i}/Z_{i})\beta+{(\tau_{i}/Z_{i})}^{1/2}\epsilon_{i}
\label{data2}\EndEq

Now using a change of variable $w_{i}=z_{i}/\tau_{i}$ results in a
likelihood of $\X|\mu,\beta,\theta$ expressible as \Eq
\Lcr(\X|\mu,\beta,\theta)=\int_{{\mathbb R}^{n}_{+}}
\[\prod_{i=1}^{n}\phi(X_{i}|\mu\Delta+\beta
w^{-1}_{i},w^{-1}_{i})\]\E\[\prod_{i=1}^{n}{\mbox
e}^{-w_{i}\tau_{i}}\tau_{i}\]\prod_{i=1}^{n}dw_{i}\label{like2}
\EndEq

where \Eq \E\[\prod_{i=1}^{n}{\mbox
e}^{-w_{i}\tau_{i}}\tau_{i}\]=\int_{{\mathbb
R}^{n}_{+}}\[\prod_{i=1}^{n}{\mbox e}^{-w_{i}\tau_{i}}\tau_{i}\]
f(\tau_{1},\ldots,\tau_{n}|\theta)\prod_{i=1}^{n}d\tau_{i}\label{james}\EndEq

It is noteworthy that~\mref{like2} also has the form, $$
\Lcr(\X|\mu,\beta,\theta)=C(\X|\mu,\beta)\int_{{\mathbb
R}^{n}_{+}} \E\[\prod_{i=1}^{n}{\mbox
e}^{-w_{i}\tau_{i}}\tau_{i}\]\prod_{i=1}^{n}f_{GIG}(w_{i}|3/2,|\beta|,|A_{i}|)dw_{i}
$$ where $C(\X|\mu,\beta)$ is determined from the Normal density
and the GIG density. As we shall show the expression
in~\mref{james} is easily handled by applying the results
James~(2005a, 2002). In closing this section notice that once
$\tau$ is integrated out in~\mref{like2} that one has a model
whereby $X_{i}|w_{i},\mu,\beta$ are independent \Eq
X_{i}=\mu\Delta+w^{-1}_{i}\beta+{w_{i}}^{-1/2}\epsilon_{i}
\label{data3}\EndEq Hence conditional on
$\W=(W_{1},\ldots,W_{n})$, the parameters $(\mu,\beta)$ are easily
estimated by standard parametric methods. \Remark Note that the
likelihood models, ~\mref{BNSlik}, analyzed in James~(2005b),
depended on $\tau$ only through terms such as
$E\[\prod_{i=1}^{n}{\mbox e}^{-w_{i}\tau_{i}}\]$ rather than
$\E\[\prod_{i=1}^{n}{\mbox e}^{-w_{i}\tau_{i}}\tau_{i}\].$ In
analogy to survival analysis, the first expression can be thought
of as the likelihood of a model where one only observes
right-censored observations, where the latter may represent the
appearance of both complete and censored observations. This
creates a fundamental difference in their respective marginal
analysis and structure.\EndRemark
\subsection{Model flexibility}
One important advantage of the present approach, over say the
direct use of~\mref{BNSlik}, is that we can more easily handle
variations in the model. Briefly, we mention the following
variation, which we believe helps address a question raised by
M.C. Jones in the discussant section of Barndorff-Nielsen and
Shephard~(2001a, p. 225, 237), $$ X_{i}\overset
{d}=\mu\Delta+{(\tau_{i}/Z_{i})}^{1/2}\beta+{(\tau_{i}/Z_{i})}^{1/2}\epsilon_{i}.
$$ Owing to the same derivations above, this leads to a model
whereby the $X_{i}|w_{i},\mu,\beta$ are independent $$
X_{i}=\mu\Delta+w^{-1/2}_{i}\beta+{w_{i}}^{-1/2}\epsilon_{i}.$$
More generally, for known real numbers $(a_{0},a_{1},\ldots,
a_{k})$ and possibly unknown $(\beta_{1},\ldots,\beta_{k})$ our
approach extends to models of the type
$$
X_{i}\overset
{d}=\mu\Delta+\sum_{j=1}^{k}{(\tau_{i}/Z_{i})}^{a_{j}}\beta_{j}+{(\tau_{i}/Z_{i})}^{a_{0}}\epsilon_{i}
$$
and beyond. This is seen by the fact the transformation
$w_{i}=\tau_{i}/z_{i}$, yields the models $X_{i}|w_{i},\mu,\beta$
are independent
$$
X_{i}=\mu\Delta+\sum_{j=1}^{k}{w_{i}}^{-a_{j}}\beta_{j}+{w_{i}}^{-a_{0}}\epsilon_{i}
$$
Note that the marginal distribution of
$\W|\mu,\beta_{1},\ldots,\beta_{k}$ is the same as for the case
of~\mref{data3}. To be quite clear, all our forthcoming results
hold for this more general setting by replacing
$\phi(X_{i}|\mu\Delta+\beta w^{-1}_{i},w^{-1}_{i})$ with
$$
\phi(X_{i}|\mu\Delta+\sum_{j=1}^{k}w^{-a_{j}}_{i}\beta_{j},{w_{i}}^{-2a_{0}}).$$
Obviously in this case the density~\mref{GHSV} has to be replaced
by a more general Normal-Gamma mixture, but is otherwise just as
easy to implement. \Remark It is not difficult to deal with models
where say $ \tau_{i}{\mbox e}^{-\tau_{i}} $ is replaced by
$\tau^{\alpha}_{i}{\mbox e}^{-\tau_{i}}$ for $0<\alpha<1$. Based
on our approach one simply writes
$$
\tau^{\alpha}_{i}=\tau_{i}\tau^{-(1-\alpha)}_{i}=\frac{\tau_{i}}{\Gamma(1-\alpha)}\int_{0}^{\infty}{\mbox
e}^{-y_{i}\tau_{i}}y_{i}^{-\alpha}dy_{i}
$$
An augmentation reveals that the likelihood would involve an
additional $n$ latent variables. \EndRemark

\section{Evaluation of the likelihood for general $\tau$}
Similar to James~(2005b) we now evaluate the likelihood in the
case where $\tau$ are more generally modeled as linear functionals
of a Poisson random measure defined over Polish spaces. Let $N$
denote a Poisson random measure on some Polish space $\Vcr$ with
mean intensity,
$$
\E[N(dv)]=\nu(dv).
$$
We denote the Poisson law of $N$ with intensity $\nu$ as
$\Pe(dN|\nu)$. The Laplace functional for $N$ is defined as
$$
\E[{\mbox e}^{-N(f)}]=\int_{\M}{\mbox e}^{-N(f)}\Pe(dN|\nu)={\mbox
e}^{-\Lambda(f)}
$$
where for any positive $f$, $N(f)=\int_{\Vcr}f(x)N(dx)$ and
$\Lambda(f)=\int_\Vcr(1-{\mbox e}^{-f(x)})\nu(dx).$ $\M$ denotes
the space of boundedly finite measures on $\Vcr$ [see Daley and
Vere-Jones~(1988)].
 We suppose that
$\tau_{i}=N(f_{i})$, for $i=1,\ldots, n$ where
$f_{1},\ldots,f_{n}$ are positive measurable functions on $\Vcr$.
Notice now that the index $i=1,\ldots,n$ need not correspond to
fixed intervals involving $\Delta$. With this in mind, let
$(w_{1},\ldots,w_{n})$ denote arbitrary non-negative numbers. We
shall assume throughout that $f_{1},\ldots,f_{n}$ are such that
$\Lambda(\sum_{i=1}^{n}w_{i}f_{i})<\infty$. Notice first that one
can write
$$
\prod_{i=1}^{n}\tau_{i}{\mbox e}^{-w_{i}\tau_{i}}={\mbox
e}^{-N(\sum_{i=1}^{n}w_{i}f_{i})}\[\prod_{i=1}^{n}\int_{\Vcr}f(V_{i})N(dV_{i})\].
$$
Removing the integrals, one can treat the
$\V=(V_{1},\ldots,V_{n})$ as missing values taking their values in
$\Vcr$. Similar to the case of the Blackwell-MacQueen
distribution, which plays a fundamental role in Dirichlet and
Gamma process mixture models~[see Lo (1984) and Lo and
Weng~(1989), Ishwaran and James~(2004) and James~(2005a)], we can
express the $\V$ as follows. Let
$\V^{*}=(V^{*}_{1},\ldots,V^{*}_{n(\p)})$ denote the $n(\p)\leq
n$, distinct values of $\V$, where $\p=\{C_{1},\ldots,C_{n(\p)}\}$
denotes a partition of the integers $\{1,2,\ldots,n\}$,with cells
$C_{j}=\{i:V_{i}=V^{*}_{j}\}$ for $j=1,\ldots,n(\p)$.
Additionally, let $e_{j,n}$, sometimes written as $e_{j}$, denote
the size, or cardinality, of the cell $C_{j}$. Define
$\Omega_{n}(v)=\sum_{i=1}^{n}w_{i}f_{i}(v)$. Let
$\Pe(dN|\nu_{\Omega_{n}},\V)$ correspond to the law of the random
measure $N_{\Omega_{n}}+\sum_{j=1}^{n(\p)}\delta_{V^{*}_{j}},$
where conditional on $(\V,\W)$, $N_{\Omega_{n}}$ is a Poisson
random measure with mean intensity $\nu_{\Omega_{n}}(dv):={\mbox
e}^{-\Omega_{n}(v)}\nu(dv).$ Now an application of James~(2005a,
Proposition 2.3), augmenting~\mref{like2}, yields a joint
distribution of $(\V,\W,N,\X)$ which is expressible as, \Eq
\Pe(dN|\nu_{\Omega_{n}},\V)\[\prod_{i=1}^{n}f(V_{i})\]\[\prod_{j=1}^{n(\p)}{\mbox
e}^{-\Omega_{n}(V^{*}_{j})}\nu(dV^{*}_{j})\]{\mbox
e}^{-\Lambda(\Omega_{n})}\prod_{i=1}^{n}\phi(X_{i}|\mu\Delta+\beta
W^{-1}_{i},W^{-1}_{i})\label{keyprop}\EndEq where
$${\mbox e}^{-\Lambda(\Omega_{n})}=\E\[{\mbox
e}^{-N(\sum_{i=1}^{n}w_{i}f_{i})}\]=\int_{\M}{\mbox
e}^{-N(\sum_{i=1}^{n}w_{i}f_{i})}\Pe(dN|\nu).$$
 Note also that
$$\Mcr(d\V|\nu_{\Omega_{n}})=\nu_{\Omega_{n}}(dV_1)\prod_{i=2}^{n}
\[\nu_{\Omega_{n}}(dV_{i})+\sum_{j=1}^{n(\p_{i-1})}\delta_{V^{*}_{j}}(dV_{i})\]=\prod_{j=1}^{n(\p)}{\mbox
e}^{-\Omega_{n}(V^{*}_{j})}\nu(dV^{*}_{j})$$ corresponds to the
n-th moment measure of a Poisson random measure with intensity
$\nu_{\Omega_{n}}$, and importantly has a structure similar to the
Blackwell-MacQueen distribution. The expression $\p_{i-1}$
corresponds to a partition of the integers $\{1,\ldots, i-1\}$ for
$i\ge 2$. Now integrating out $(\V,\W,N)$ in~\mref{keyprop} leads
to an expression for the likelihood.

\begin{thm} Suppose that $\tau_{i}=N(f_{i})$ for $i=1,\ldots,n$
where $N$ is a Poisson random measure on $\Vcr$ with intensity
$\nu$. Then the likelihood~\mref{like2} can be expressed as
$$\Lcr(\X|\mu,\beta,\theta)=\int_{{\mathbb R}^{n}_{+}}\Bcr({\mathbf
w}){\mbox
e}^{-\Lambda(\Omega_{n})}\prod_{i=1}^{n}\phi(X_{i}|\mu\Delta+\beta
w^{-1}_{i},w^{-1}_{i})dw_{i}
$$
where $\Bcr({\mathbf
w})=\int_{{\Vcr}^{n}}\prod_{i=1}^{n}f_{i}(V_{i})\Mcr(d\V|\nu_{\Omega_{n}})=\sum_{\p}
\prod_{j=1}^{n(\p)}\int_{\Vcr}\[\prod_{i\in
C_{j}}f_{i}(v)\]\nu_{\Omega_{n}}(dv)$. Where $\sum_{\p}$ denotes
the sum over all partitions of the integers $\{1,\ldots,
n\}.$\qed\end{thm}

\begin{thm} Suppose that $\tau_{i}=N(f_{i})$ for $i=1,\ldots,n$
where $N$ is a Poisson random measure on $\Vcr$ with intensity
$\nu$. Then augmenting the likelihood in Theorem 3.1 yields the
joint posterior distribution of $\V,\W|\X$ given by,
$$
\Pcr_{\X}(d\v,d\w)\propto\[\prod_{i=1}^{n}f_{i}(v_{i})\]\Mcr(d\V|\nu_{\Omega_{n}}){\mbox
e}^{-\Lambda(\Omega_{n})}\prod_{i=1}^{n}\phi(X_{i}|\mu\Delta+\beta
w^{-1}_{i},w^{-1}_{i})dw_{i}
$$
In particular we have the following posterior distributions
\Enumerate
\item[(i)]$\pi(d\w|\v,\X)\propto {\mbox e}^{-\Lambda(\Omega_{n})}\prod_{i=1}^{n}
f_{GIG}(w_{i}|3/2,|\beta|,\sqrt{A^{2}_{i}+2\sum_{j=1}^{n(\p)}f_{i}(v^{*}_{j})})$
\item[(ii)]$\pi(d\v|\w,\X)\propto \[\prod_{i=1}^{n}f(v_{i})\]\Mcr(d\v|\nu_{\Omega_{n}})=\[\prod_{j=1}^{n(\p)}
\prod_{i\in C_{j}}f_{i}(v^{*}_{j})\]\Mcr(d\v|\nu_{\Omega_{n}})$
\item[(iii)]$\Pcr_{\X}(d\w)=\Bcr({\mathbf w}){\mbox
e}^{-\Lambda(\Omega_{n})}\prod_{i=1}^{n}\phi(X_{i}|\mu\Delta+\beta
w^{-1}_{i},w^{-1}_{i})dw_{i}/\Lcr(\X|\mu,\beta,\theta)$ is the
posterior density of $\W|\X.$ \EndEnumerate \qed
\end{thm}

\begin{prop}The distribution of $\V|\W,\X$ can be further described
as follows. The posterior distribution of $\V|\p,\W,\X$, is such
that the unique values $\V^{*}$ are conditionally independent with
distributions
$$
\Pe(V^{*}_{j}\in dv|\p,\w,\X)\propto\[\prod_{i\in
C_{j}}f_{i}(v)\]\nu_{\Omega_{n}}(dv){\mbox { for
}}j=1,\ldots,n(\p).
$$
The posterior distribution of $\p|\W,\X$ is given by
$\pi(\p|\w,\X)\propto\prod_{j=1}^{n(\p)}\int_{\Vcr}\[\prod_{i\in
C_{j}}f_{i}(v)\]\nu_{\Omega_{n}}(dv).$\qed
\end{prop}

\subsection{Posterior distribution of parameters}
Upon examining the likelihood, one sees that Bayesian inference
for parameters $(\mu,\beta,\theta)$ can be implemented in a
straightforward manner, along the lines of methods outlined for
the Dirichlet/Gamma process semi-parametric mixture models. See
Ishwaran and James~(2004) for these ideas and further pertinent
references. Specifically, a straightforward application of Bayes
rule yields the following results.
\begin{prop} Suppose that $\tau$ depends on a $d$-dimensional
parameter $\theta$. Then if $q(d\theta)$, $q(d\beta)$, $q(\mu)$
denote independent priors for $(\beta,\mu,\theta)$, their
posterior distributions can be written as follows \Enumerate
\item[(i)]$q(d\beta|\mu,\w,\X)\propto q(d\beta){\mbox
e}^{-[\beta^{2}/2\sum_{i=1}^{n}w_{i}-n{\bar A}\beta]}$
\item[(ii)]$q(d\mu|\beta,\w,\X)\propto q(d\mu){\mbox
e}^{-[\sum_{i=1}^{n}A^{2}_{i}w^{-1}_{i}-n{\bar A}\beta]}$
\item[(iii)]$q(d\theta|\w,\V,\X)\propto q(d\theta){\mbox
e}^{-\Lambda_{\theta}(\Omega_{n})}\prod_{j=1}^{n(\p)}\nu_{\theta}(dV^{*}_{j})$
\EndEnumerate\qed
\end{prop}
\subsection{Posterior distribution of the process}
The above results describe the behavior of the finite-dimensional
likelihood and parameters. It is useful to also obtain a
description of the underlying random process given the data. This
allows one to see directly how the data affects the overall
process. Moreover, combined with the results in James~(2005a), it
provides a calculus for more general functionals. For notational
simplicity we suppose that $(\mu,\beta,\theta)$ are fixed. The
next result also follows immediately from an application of
Fubini's theorem and~\mref{keyprop}.
\begin{thm} Suppose that a likelihood of $\X$ and the specifications for $\tau$ and $N$ are defined by the
specifications in Theorem 3.1. Let $\Omega_{n}(x)=\sum_{i=1}W_{i}
f_{i}(x)$. Let $\Scr=\Vcr\times(0,\infty)$. Then the posterior
distribution of $N|\X$ is
$$\int_{\Scr^{n}} \Pe(dN|\nu_{\Omega_{n}},\v)\Pcr_{\X}(d\v,d\w).$$
In particular for any positive or integrable function $g$ on $\M$,
$$\int_{\Scr^{n}}\[\int_{\M}g(N)\Pe(dN|\nu_{\Omega_{n}},\v)\]\Pcr_{\X}(d\v,d\w)
=\int_{\Scr^{n}}\[\int_{\M}g(N+\sum_{j=1}^{n(\p)}\delta_{v^{*}_{j}})\Pe(dN|\nu_{\Omega_{n}})\]
\Pcr_{\X}(d\v,d\w)
$$
\qed
\end{thm}

\subsection{A general posterior predictive density for the log price}
We now define a random variable similar to~\mref{genrv} which can
be thought of as representing the log-price and give an explicit
expression for its posterior density given $\X$. The random
variable is defined as, $${\tilde X}\overset {d}=\mu{\tilde
\Delta}+({\tilde \tau}/Z)\beta+{({\tilde \tau}/Z)}^{1/2}{\tilde
\epsilon},$$ where ${\tilde \Delta}$ is just some non-negative
number, ${\tilde \tau}=N({\tilde f})$ for some positive function
${\tilde f}$ such that Laplace transform of ${\tilde \tau}$
exists. $\epsilon$ is a standard Normal random variable. Now let
$\Omega_{n+1}(x)=\Omega_{n}(x)+w{\tilde f}(x)$.

\begin{prop}The posterior density of the log stock price given
${\mathbf X}$ is, $f_{{\tilde X}}(x|\beta,\mu,{\mathbf X})$ equal
to,
$$
\int_{{\mathbb R}^{n}_{+}}\left[\int_{0}^{\infty}q(w|\v,\w){\mbox
e}^{-[\Lambda(\Omega_{n+1})-\Lambda(\Omega_{n})]}
\phi(x|\mu{\tilde \Delta}+\beta
w^{-1},w^{-1})dw\right]\Pcr_{{\mathbf X}}(d\v,d{\mathbf w}),$$
where $q(w|\V,\w)=\[\int_{\Vcr}{\tilde f}(v){\mbox e}^{-w{\tilde
f}(v)}\nu_{\Omega_{n}}(dv)+\sum_{j=1}^{n(\p)}{\tilde
f}(V^{*}_{j})\]{\mbox e}^{-\sum_{j=1}^{n(\p)}w{\tilde
f}(V^{*}_{j})}.$\qed
\end{prop}
\Proof The result follows from Theorem 3.3, and~\mref{mixrep},
using the following fact,
$$
\int_{\M}N({\tilde f}){\mbox e}^{-N(w{\tilde
f})}\Pe(dN|\nu_{\Omega_{n}},\V)={\mbox
e}^{-\sum_{j=1}^{n(\p)}w{\tilde f}(V^{*}_{j})}\int_{\M}[N({\tilde
f})+\sum_{j=1}^{n(\p)}{\tilde f}(V^{*}_{j})]{\mbox e}^{-N(w{\tilde
f})}\Pe(dN|\nu_{\Omega_{n}}),
$$
where $wN({\tilde f})=w{\tilde \tau}$, and $\int_{\M}{\mbox
e}^{-N(w{\tilde f})}\Pe(dN|\nu_{\Omega_{n}})={\mbox
e}^{-[\Lambda(\Omega_{n+1})-\Lambda(\Omega_{n})]}$. \EndProof

\subsection{Generalized Chinese Restaurant and P\'olya Urn procedures}
The results in the previous sections show that, viewing $\W$ as a
parameter, these models are structurally similar to Bayesian
semi-parametric mixture models based on multiplicative intensity
likelihoods. As such, one can import computational sampling
schemes described in Ishwaran and James~(2004) and James~(2005a)
and references therein. One can deduce the necessary modifications
for the  general nonparametric setting in James~(2005a)~from the
methods for the semi-parametric Gamma process setting described in
Ishwaran and James~(2004).  In particular, this includes  general
semi-parametric analogues of P\'olya Urn Gibbs samplers and SIS
procedures given by Escobar~(1994), Liu~(1996), and West, M\"uller
and Escobar~(1994), and the Gibbs sampling/SIS procedures based on
a generalized weighted Chinese restaurant process~[see Lo, Brunner
and Chan~(1996) and Ishwaran and James~(2003)]. However, since
these schemes are phrased primarily for completely random
measures, which are a subclass of the models we look at here, we
mention a few details for clarification in the general Poisson
case. First note that sampling from the distribution of
$\W,\theta,\beta,\mu|\V,\X$, described by Theorem 3.2 and
Proposition 3.2,  proceeds along the lines of well-known
parametric procedures such as random walk Metropolis-Hastings. The
task then remains to approximately sample
$\V|\W,\theta,\beta,\mu$. The key fact, is that structurally these
models are not markedly different than the Dirichlet/Gamma process
mixture models based on the Blackwell-MacQueen P\'olya Urn
distribution. In fact, a weighted Chinese restaurant SIS algorithm
to sample Urn distributions derived from general Poisson random
measures, such as that of $\V$, has already been given in
James~(2002, section 2.3). This of course translates into dual
Gibbs sampling procedures. Here we sketch out the relevant
probabilities to implement these type of schemes for the
simulation from the distribution,
$$\pi(d\V|\W,\beta,\theta,\mu,\X)\propto
\[\prod_{i=1}^{n}f_{i}(V_{i})\]\Mcr(d\V|\nu_{\Omega_{n}})=\[\prod_{j=1}^{n(\p)}
\prod_{i\in C_{j}}f_{i}(V^{*}_{j})\]\Mcr(d\V|\nu_{\Omega_{n}})$$
Note that this distribution will typically not depend on $(\beta,\
\mu)$ except through $\W$. Similar to James~(2005a, equation 40),
define for $r=0,\ldots, n-1$ conditional probabilities,
    $$
    \Pe(V_{r+1}\in
    dx|\V_{r})=\frac{l_{0,r}}{c_{r}}\lambda_{r}(dx)+\sum_{j=1}^{n(\p_{r})}
    \frac{l_{j,r}(V^{*}_{j})}{c_{r}}\delta_{V^{*}_{j}}(dx)
    $$
where $\V_{r}=\{V_{1},\ldots,V_{r}\}$, $\lambda_{r}(dx)\propto
f_{r+1}(x)\nu_{\Omega_{n}}(dx)$ and
$l_{0,r}=\int_{\Vcr}f_{r+1}(x)\,\nu_{\Omega_{n}}(dx)$ and
$l_{j,r}(x)=f_{r+1}(x)$ Additionally
$c_{r}=l_{0,r}+\sum_{j=1}^{n(\p_{r})}l_{j,r}(V^{*}_{j})$.
Examining James~(2005a, section 4.4.) we see these are the
ingredients to implement general analogues of the P\'olya Urn
Gibbs Sampler and SIS procedures described by Escobar~(1994) and
Liu~(1996). An acceleration step similar to West, M\"uller and
Escobar~(1994) can be implemented by using the description in
Ishwaran and James~(2004, p.180, Remark 2) combined with
Proposition 3.1. Naturally, if one has a structure closer to
completely random measures the P\'olya Urn type methods described
in James~(2005a, section 4.4.), which involve integrating the jump
components in the $\V$ vector should be employed if possible. To
get the Chinese restaurant type procedures one samples partitions
$\p$ based on probabilities derived from $l_{0,r}$ and
$$
l_{j,r}=\int_{\Vcr}l_{j,r}(x)\,
    \[\prod_{i\in C_{j,r}}f_{i}(x)\]\nu_{\Omega_{n}}(dx){\mbox { for }}
    j=1,\ldots,n(\p_{r})
$$
where $\p_{r}$ denotes a partition of the integers $\{1,\ldots,
r\}$ and each $C_{j,r}=\{i\leq r: V_{i}=V^{*}_{j}\}$ denotes the
corresponding cells. Additionally,
$l(r)=l_{0,r}+\sum_{j=1}^{n(\p_{r})}l_{j,r}$. See James~(2002,
section 2.3) for justification of these procedures. \Remark The
main distinction, structurally, between the nonparametric
multiplicative intensity models described in James~(2005a) and the
present {\it semi-parametric} setting, is the non-cancelation of
the Laplace transform ${\mbox e}^{-\Lambda(\Omega_{n})}$ in the
likelihood as it depends on external parameters $(\W,\theta)$.
This point is addressed in Ishwaran and James~(2004) for Gamma
processes. \EndRemark

\section{Results for completely random measures} Many processes
$\tau$ will be directly expressible as functionals of completely
random measures, say $\mu$. This is the case for models based on
the integrated OU processes of Barndorff-Nielsen and
Shephard~(2001a,b), where $\mu$ is the {\it Background Driving
L\'evy Process}~(BDLP). As such, refinements of the above results
in that case can be deduced from James~(2005a, section 4). Note
that a homogeneous completely random measure, say $\mu$, with no
drift, has the representation
$\mu(dy)=\int_{0}^{\infty}uN(du,dy)$. Where,
$\Vcr=(0,\infty)\times\Ycr$ for some Polish space $\Ycr$,
$\nu(du,dy):=\rho(du)\eta(dy)$. The measure $\rho$ is the L\'evy
density of a non-negative infinitely divisible random variable,
$T$, with Laplace transform
$$
\E\[{\mbox e}^{-\omega T}\]={\mbox e}^{-\psi(\omega)}
$$
where $\psi(\omega)=\int_{0}^{\infty}(1-{\mbox e}^{-\omega
u})\rho(du)$. One may then write $V_{i}=(J_{i},Y_{i})$ and
$V^{*}_{j}=(J_{j,n},Y^{*}_{j})$, where $(J_{j,n})$ denotes the
unique jumps of the process $\mu$, picked by a type of biased
sampling. Moreover it follows that for measureable functions
$f_{i}(u,y)=ug_{i}(y)$, $\tau_{i}:=N(f_{i})=\mu(g_{i})$.
Additionally for any $f$ and $g$ such that $f(u,y)=ug(y)$, one has
$\Lambda(f)=\int_{\Ycr}\psi(g(y))\eta(dy)<\infty$. Define,
$\rho_{\Omega_{n}}(du|y)={\mbox
e}^{-u\sum_{i=1}^{n}w_{i}g_{i}(y)}\rho(du)$, and for
$l=1,\ldots,n$ {\it conditional} cumulants \Eq
\kappa_{l}(\rho_{\Omega_{n}}|y)=\int_{0}^{\infty}u^{l}\rho_{\Omega_{n}}(du|y).
\label{cumulants} \EndEq

\begin{thm} Suppose that $N$ is a Poisson random measure with
intensity $\nu(du,dy)=\rho(du)\eta(dy)$
on~$\Vcr=(0,\infty)\times\Ycr$. Suppose that
$\tau_{i}:=N(f_{i})=\mu(g_{i})$ for $i=1,\ldots,n$. Then according
to the model~\mref{like2}, one has the following results.
\Enumerate
\item[(i)] Setting $V^{*}_{j}=(J_{j,n},Y^{*}_{j})$, for
$j=1,\ldots, n(\p)$, conditional on $\p,\W,\X$, the pairs of
random variables on $\Vcr$ are independent with distributions
$$\Pe(J_{j,n}\in du,Y^{*}_{j}\in dy|\w,\X)\propto
\frac{u^{e_{j}}\rho_{\Omega_{n}}(du|y)}
{\kappa_{e_{j}}(\rho_{\Omega_{n}}|y)}\kappa_{e_{j}}(\rho_{\Omega_{n}}|y)\[\prod_{i\in
C_{j}}g_{i}(y)\]\eta(dy),$$ where $ \int_{\Vcr}\[\prod_{i\in
C_{j}}f_{i}(v)\]\nu_{\Omega_{n}}(dv)=\int_{\Vcr}\kappa_{e_{j}}(\rho_{\Omega_{n}}|y)\[\prod_{i\in
C_{j}}g_{i}(y)\]\eta(dy):=\vartheta(C_{j}|\w)$ is the normalizing
constant.
\item[(ii)]The posterior distribution of $\mu|\p,\W,\X$, is such that
$\mu(dx)\overset
{d}=\mu_{\Omega_{n}}(dx)+\sum_{j=1}^{n(\p)}J_{j,n}\delta_{Y^{*}_{j}}(dx)$
where given $\p,\W,\X$, $\mu_{\Omega_{n}}$ is a completely random
measure determined by the law $\Pe(dN|\nu_{\Omega_{n}})$, and the
pairs $(J_{j,n},Y^{*}_{j})$ are conditionally independent of
$\mu_{\Omega_{n}}$ with distribution described in $[(i)]$
\item[(iii)]The distribution of $\p|\W,\X$ is proportional to
$\[\prod_{j=1}^{n(\p)}{\vartheta(C_{j}|\w)}\].$
\item[(iv)]The density of $\W|\p,\X$ is $$f(\w|\p,\X)\propto{\mbox e}^{-\int_{\Ycr}\psi(\sum_{i=1}^{n}w_{i}g_{i}(y))\eta(dy)}\[\prod_{j=1}^{n(\p)}
{\vartheta(C_{j}|\w)}\]\prod_{i=1}^{n}\phi(X_{i}|\mu\Delta+\beta
w^{-1}_{i},w^{-1}_{i})$$ \EndEnumerate\ All the above results
hold, in an obvious way, for the inhomogeneous case of
$\nu(du,dy)=\rho(du|y)\eta(dy).$\qed
\end{thm}
\subsection{Remarks on implementation} The following remarks
address specifically the models in this section, but clearly have
extensions to the general setting. Suppose the infinitely
divisible random variable $T$ has density $f_{T}$, and hence for a
unique $\rho$, it has Laplace transform $\int_{0}^{\infty}{\mbox
e}^{-\omega t}f_{T}(t)dt={\mbox e}^{-\psi(\omega)}$. It is then
important to note that $\kappa_{l}(\rho_{\Omega_{n}}|y)$ are for
fixed $y$, the first $l=1,\ldots,n$ cumulants of an infinitely
divisible random variable with density,
$$f_{\Omega_{n}}(t|y):={\mbox
e}^{-t\sum_{i=1}^{n}w_{i}g_{i}(y)}f_{T}(t){\mbox
e}^{\psi(\sum_{i=1}^{n}w_{i}g_{i}(y))}.$$ Now defining the the
first $l=1,\ldots,n$ moments as
$m_{l}(y)=\int_{0}^{\infty}t^{l}f_{\Omega_{n}}(t|y)dv$, it follows
that the cumulants may be calculated using the result of Theile.
That is,  $$
    \kappa_{l}(\rho_{
    \Omega_{n}}|y)=m_{l}(y)-\sum_{k=1}^{l-1}\binom{l-1}{k-1}
    \kappa_{k}(\rho_{
    \Omega_{n}}|y)m_{l-k}(y).
    $$
 This indicates quite
clearly, the important fact, that one need not have the specific
form of the L\'evy density $\rho$ to implement estimation
procedures for our models. An interesting case would be where $T$
has a log Normal distribution. Of course for models such as stable
laws where $\rho$ has a simple form and the probability density is
generally complex, the converse is also true.
\subsection{Example: Generalized Gamma processes}
We now provide some details for one of the most tractable classes
of models. An interesting  class of measures are the family of
generalized Gamma random measures discussed in Brix~(1999). Using
the description of Brix~(1999), these are $\mu$ processes with
L{\'e}vy measure
$$\rho(du)=\frac{1}{
\G(1-\a)}u^{-\alpha-1}{\mbox e}^{-bu}du.
$$
The values for $\alpha$ and $b$ are restricted to satisfy $0<\a<1$
and $0\le b<\infty$ or $-\infty<\a\le 0$ and $0<b<\infty$.
Different choices for $\a$ and $b$ yield various subordinators.
These include the stable subordinator when $b=0$, the Gamma
process subordinator when $\alpha=0$ and the inverse-Gaussian
subordinator when $\a=1/2$ and $b>0$. When $\alpha <0$ this
results in a class of Gamma compound Poisson processes. It follows
that, for $\alpha\neq 0$, and $b\ge 0$, $
\psi(\sum_{i=1}^{n}w_{i}g_{i}(y))={\frac{1}{\alpha}}[(b+\sum_{i=1}^{n}w_{i}g_{i}(y))^{\alpha}-{b}^{\alpha}]
$ The case of the Gamma process, $\alpha=0$, $b>0$ is a limiting
case and results in the well-known expression $
\psi(\sum_{i=1}^{n}w_{i}g_{i}(y))=\ln(1+\sum_{i=1}^{n}w_{i}g_{i}(y)/b)
$ Now, for all $\alpha$ and $b$, conditional on
$Y^{*}_{j},\p,\W,\X$, each $J_{j,n}$ is Gamma distributed with
shape and scale parameters
$(e_{j}-\alpha,b+\sum_{i=1}^{n}w_{i}g_{i}(Y^{*}_{j}))$. It follows
that the joint moment measure of $\Y|\W,\X$ is, $$
\[\prod_{j=1}^{n(\p)}
\frac{\Gamma(e_{j,n}-\alpha)}{\G(1-\a)}\]\prod_{j=1}^{n(\bf
  p)}{(b+\sum_{i=1}^{n}w_{i}g_{i}(Y^{*}_{j}))}^{-(e_{j,n}-\alpha)}\eta(dY^{*}_{j})
$$
which is the key component in the sampling algorithms described in
James~(2005a). See Ishwaran and James~(2004) for many details,
including the usage of {\it Blocked Gibbs Samplers}, in the Gamma
process semi-parametric setting which easily translates to the
generalized Gamma class.
\section{Example: BNS-OU model}
For proper comparison with the models~\mref{BNSlik} as derived in
James~(2005b), it is interesting to look at how the Integrated OU
model of Barndorff-Nielsen and Shephard~(2001a,b, 2003) behave in
this scenario.  We shall refer to this model as the BNS-OU model.
Throughout this section we shall take $\Ycr=(-\infty,\infty)$. One
may express the Barndorff-Nielsen and Shephard~(2001 a, b)
integrated OU process $\tau$ as \Eq \tau(t)=\lambda^{-1}[(1-{\mbox
e}^{-\lambda t})\int_{-\infty}^{0}{\mbox
e}^{y}\mu(dy)+\int_{0}^{t}(1-{\mbox
e}^{-\lambda(t-y)})\mu(dy)]\label{model1} \EndEq where
$v(0):=v_{0}=\int_{-\infty}^{0}{\mbox e}^{y}\mu(dy)$, denotes the
{\it instantaneous volatility} at time $0.$ The form in
~\mref{model1} is taken from Carr, Geman, Madan and Yor~(2003, p.
365). It follows that for any $s<t$,
$[\tau(t)-\tau(s)]=\mu(g_{s,t})=N(f_{s,t})$ where
$f_{s,t}(u,y)=ug_{s,t}(y)$ and $\lambda g_{s,t}(y)$ equals,

\Eq {\mbox e}^{-\lambda s}(1-{\mbox e}^{-\lambda (t-s)}){\mbox
e}^{y}I_{\{y\leq 0\}}+(1-{\mbox e}^{-\lambda(t-y)})I_{\{s<y \leq
t\}}+{\mbox e}^{-\lambda s}(1-{\mbox e}^{-\lambda(t-s)}){\mbox
e}^{\lambda y}I_{\{0<y\leq s\}}.\label{gendif3}\EndEq

The first component in~\mref{gendif3} represents the contribution
from $v_{0}$. Specializing this to $s=(i-1)\Delta$ and $t=i\Delta$
one has $\tau_{i}=\mu(g_{i})=N(f_{i})$ where
$f_{i}(u,y)=ug_{i}(y)$ and further
$g_{i}(y)=g_{i,1}(y)+g_{i,2}(y)$ with, \Eq
g_{i,1}(y)=\lambda^{-1}[(1-{\mbox e}^{-\lambda
(i\Delta-y)})I_{\{(i-1)\Delta<y\leq i\Delta\}}+{\mbox
e}^{-\lambda(i-1)\Delta}(1-{\mbox e}^{-\lambda\Delta}){\mbox
e}^{y}I_{\{y\leq 0\}}]\label{gform} \EndEq and \Eq
g_{i,2}(y)=\lambda^{-1}{\mbox e}^{-\lambda(i-1)\Delta}(1-{\mbox
e}^{-\lambda \Delta}){\mbox e}^{\lambda y}I_{\{0<y\leq
(i-1)\Delta\}}.\label{hform}\EndEq Now for $i=1,\ldots, n$, set
$r_{i}=\lambda^{-1}[\sum_{k=i}^{n}w_{k}{\mbox e}^{-\lambda (k-1)
\Delta}](1-{\mbox e}^{- \lambda \Delta})$, and define $r_{n+1}=0$.

Now notice that for any sequence of numbers, the simplest
expression will be obtained by utilizing the following facts. \Eq
\sum_{j=1}^{n}w_{j}[g_{j,1}{(y)}+g_{j,2}(y)]=r_{1}{\mbox
e}^{y}{\mbox { for }}y\leq 0 \label{id1}\EndEq and for
$i=1,\ldots, n$\Eq
\sum_{j=1}^{n}w_{j}[g_{j,1}{(y)}+g_{j,2}(y)]=\zeta
(y|w_{i},r_{i+1}){\mbox { for }}(i-1)\Delta<y\leq i\Delta.
\label{id2}\EndEq Where for each $i$, $\zeta
(y|w_{i},r_{i+1})=[\lambda^{-1}w_{i}(1-{\mbox
e}^{-\lambda(i\Delta-y)})+r_{i+1}{\mbox e}^{\lambda y}].$ Then one
has the following result which is a generalized version of
James~(2005b, Proposition 3.1)

\begin{prop} For $0\leq s<t$, let $\tau(t)-\tau(s)$ be defined by~\mref{model1} and \mref{gendif3}
Then the results of Theorem 4.1 hold with
$f_{i}(u,y)=u[g_{i,1}(y)+g_{i,2}(y)]$, as described
in~\mref{gform} and~\mref{hform}. Suppose that
$\eta(dy):=\eta(y)dy$, and define
$\eta_{\lambda}(i\Delta,u)=\eta(i\Delta+\lambda^{-1}\ln(1-u))$.
Now, in particular, using a change of variable, \Enumerate
\item[(i)]
$ {\mbox e}^{-\Lambda(\sum_{i=1}^{n}w_{i}f_{i})}={\mbox
e}^{-\Phi_{0}(r_{1})}{\mbox
e}^{-\Phi_{n}(w_{n})}\prod_{i=1}^{n-1}{\mbox
e}^{-\Phi_{i}(w_i|r_{i+1})} $ \item[(ii)]
$\Phi(w_{i}|r_{i+1})=\int_{1-{\mbox
e}^{-\lambda\Delta}}^{1}\lambda^{-1}\psi(r_{i+1}{\mbox e}^{\lambda
i\Delta}(1-u)+\lambda^{-1}w_{i}u)\frac{\eta_{\lambda}(i\Delta,u)du}{1-u},$
for $i=1,\ldots, n-1$
\item[(iii)]
$\Phi(w_{n})=\int_{1-{\mbox
e}^{-\lambda\Delta}}^{1}\lambda^{-1}\psi(\lambda^{-1}w_{n}u)\frac{\eta_{\lambda}(i\Delta,u)du}{1-u}$
\item[(iv)]
$\Phi_{0}(r_{1})=\int_{0}^{1}\psi(r_{1}u)\eta(\ln u)\frac{du}{u}$,
where ${\mbox e}^{-\Phi_{0}(r_{1})}=\E[{\mbox e}^{-r_{1}v_{0}}].$
\EndEnumerate \qed
\end{prop}

Additionally one has the following features which do not play a
role in the analysis of James~(2005b). First statements~\mref{id1}
and ~\mref{id2} imply that $$\rho_{\Omega_{n}}(du|y)=\[{\mbox
e}^{-ur_{1}{\mbox e}^{y}}I_{\{y\leq 0\}}+\sum_{i=1}^{n}{\mbox
e}^{-u\zeta (y|w_{i},r_{i+1})} I_{\{(i-1)\Delta<y\leq
i\Delta\}}\]\rho(du).$$ Now with some abuse of notation write
$$\rho_{r_{1}}(du|y)={\mbox e}^{-ur_{1}{\mbox e}^{y}}\rho(du){\mbox { and
}}\rho_{w_{i},r_{i+1}}(du|y)={\mbox e}^{-u\zeta
(y|w_{i},r_{i+1})}\rho(du).$$ Hence for $l=1,\ldots,n$, one can
write in an obvious way,
$$
\kappa_{l}(\rho_{\Omega_{n}}|y):=\kappa_{l}(\rho_{r_{1}}|y)I_{\{y\leq
0\}}+\sum_{i=1}^{n}
\kappa_{l}(\rho_{w_{i},r_{i+1}}|y)I_{\{(i-1)\Delta<y\leq
i\Delta\}}.
$$

Let $i^{*}_{j}=\min_{i\in C_{j}}$, that is the minimal index in a
cell $C_{j}$, then,
$$
\prod_{i\in C_{j}}\left(  g_{i,1}\left(  y\right)  +g_{i,2}\left(
y\right) \right)
 =a(j,\Delta)\left(q(y|i^{*}_{j},e_{j})I_{\left\{  \left(
i_{j}^{\ast}-1\right) \Delta<y\leq i_{j}^{\ast}\Delta\right\}
}+e^{e_{j}y}I_{\left\{ y\leq0\right\} }+e^{\lambda
e_{j}y}I_{\left\{ 0<y\leq \left( i_{j}^{\ast}-1\right)
\Delta\right\}  }\right)$$ where
$$q(y|i^{*}_{j},e_{j})=\frac{(1-e^{\lambda (i_{j}^{\ast}\Delta -y)
}){\mbox e}^{\lambda (e_{j}-1)y}} {{\mbox
e}^{-\lambda(i^{*}_{j}-1)\Delta}\left( 1-e^{-\lambda\Delta}\right)
}{\mbox { and }}  a(j,\Delta)=  \prod_{i\in C_{j}}
\frac{1}{\lambda}e^{-\lambda( i-1) \Delta}(
1-e^{-\lambda\Delta}).$$ Now one has
$\kappa_{e_{j}}(\rho_{\Omega_{n}}|y)\prod_{i\in C_{j}}\left(
g_{i,1}\left(  y\right)  +g_{i,2}\left(
y\right)\right)/a(j,\Delta)$ is equal to
\begin{align*}
& r(y|i^{*}_{j},e_{j},\w)={\mbox
e}^{e_{j}y}\kappa_{e_{j}}(\rho_{r_{1}}|y)I_{\{y\leq
0\}}+q(y|i^{*}_{j},e_{j})\kappa_{e_{j}}(\rho_{w_{i^{*}_{j}},r_{i^{*}_{j}+1}}|y)I_{\{
(i_{j}^{\ast}-1)\Delta<y\leq i_{j}^{\ast}\Delta\}}\\
& +{\mbox e}^{\lambda
e_{j}y}\sum_{k=1}^{i^{*}_{j}-1}\kappa_{e_{j}}(\rho_{w_{k},r_{k+1}}|y)I_{\{(k-1)\Delta<y\leq
k\Delta\}}.
\end{align*}
One can check these points easily by looking at
$[g_{i,1}(y)+g_{i,2}(y)][g_{l,1}(y)+g_{l,2}(y)]$ for any pair
$i<l$. Additionally,
$$a_{n}:=\prod_{j=1}^{n(\p)}a(j,\Delta)=\lambda^{-n}{(1-{\mbox
e}^{\lambda \Delta})}^{n}{\mbox e}^{-\lambda n(n-1)\Delta/2}.$$
The above derivations yield the necessary specifications for the
characterization of the pertinent features of this model via
Theorem 4.1. However, our derivations above importantly reveal a
much more refined structure which has been associated with
Bayesian mixture models for monotone hazards and densities. That
is the works of Dykstra and Laud~(1981), Lo and Weng~(1989),
Brunner and Lo~(1989), Ho~(2002) and Ho~(2005).
\subsection{ A combinatorial reduction in terms of s-paths}
Similar to Brunner and Lo~(1989, p.1553), let
$\m=(m_{1},\ldots,m_{n})$ denote a vector of non-negative integers
taking values in $\Xi:=\{\m|\sum_{i=1}^{j}m_{i}\ge j, 1\leq j\leq
n-1, \sum_{i=1}^{n}m_{i}=n\}$. Then each $m_{i}$ denotes the size
of the cell whose minimal index is $i$. We also define binary
random variables $\{\xi_{1},\ldots,\xi_{n}\}$ where  $\xi_{1}:=1$
and in general $\xi_{i}=1$ if $i$ is the minimal index of a cell,
otherwise it is $0$. Viewed in terms of sequentially sampling the
random variables $\Y$, the $\xi_{i}=1$ if $Y_{i}$ is distinct from
the previous $Y_{1},\ldots,Y_{i-1}$ random variables. This idea is
a generalization of the Bernoulli random variables associated with
the Blackwell-MacQueen distribution as discussed in Korwar and
Hollander~(1973). Hence, from the derivations in the previous
section, one has
$$
\prod_{j=1}^{n(\p)}\vartheta(C_{j}|\w)=a_{n}\prod_{j=1}^{n(\p)}\int_{\Ycr}r(y|i^{*}_{j},e_{j},\w)\eta(dy)
=a_{n}\prod_{i=1}^{n}{\varphi{(i,m_{i}|\w)}}$$ where
$\varphi(i,m_{i}):=\int_{\Ycr}r(y|i,m_{i},\w)\eta(dy)$ if
$m_{i}>0$ and otherwise, $\varphi(i,0):=1$. Note importantly that
the size of the space $\Xi$ is considerably smaller than the space
of all partitions $\p$ of the integers $\{1,\ldots,n\}$.
Specifically a vector $\m$ contains information about the number
of unique values or non-empty cells
$n(\p):=\sum_{i=1}^{n}\xi_{i}$, the size of each cell $m_{i}$, and
the minimal index of each cell. However one does not know
precisely the indices $l\in C_{j}$ for each $j=1,\ldots,n(\p)$.
This leads to a more simplified version of the results in Theorem
4.1.
\begin{prop}  Suppose that $N$ is a Poisson random measure with
intensity $\nu(du,dy)=\rho(du)\eta(dy)$
on~$\Vcr=(0,\infty)\times\Ycr$. Suppose that $\tau$ is defined by
~\mref{model1}. Then the likelihood~\mref{like2} can be expressed
as
$$\Lcr(\X|\mu,\beta,\theta)=a_{n}\int_{{\mathbb R}^{n}_{+}}\[\sum_{\m\in
\Xi}\prod_{i=1}^{n}\varphi(i,m_{i}|\w)\]{\mbox
e}^{-\Phi_{0}(r_{1})}\prod_{i=1}^{n}{\mbox
e}^{-\Phi_{i}(w_i|r_{i+1})}\phi(X_{i}|\mu\Delta+\beta
w^{-1}_{i},w^{-1}_{i})dw_{i}
$$
where $\Phi_{j}$ for $j=0,\ldots,n$ are given in Proposition
5.1.\qed\end{prop}
\begin{prop} Suppose that $N$ is a Poisson random measure with
intensity $\nu(du,dy)=\rho(du)\eta(dy)$
on~$\Vcr=(0,\infty)\times\Ycr$. Suppose that $\tau$ is defined by
~\mref{model1}. Let $\MM=(M_{1},\ldots,M_{n})$ denote the random
vector corresponding to the observations $\m$. Then, one has the
following results. \Enumerate
\item[(i)] The posterior distribution of $(\J,\Y)|\MM,\W,\X$ consists of $n(\p)=\sum_{i=1}^{n}\xi_{i}$
unique values $\{({\tilde J_{i}},{\tilde Y_{i}}):\xi_i=1, 1\leq
i\leq n\}$ which are conditionally independent with respective
distributions,$$\Pe({\tilde J_{i}}\in du,{\tilde Y}_{i}\in
dy|\m,\w,\X)= \frac{u^{m_{i}}\rho_{\Omega_{n}}(du|y)}
{\kappa_{m_{i}}(\rho_{\Omega_{n}}|y)}\frac{r(y|i,m_{i},\w)\eta(dy)}{\varphi(i,m_{i}|\w)}.$$
\item[(ii)]The distribution of $\MM|\W,\X$ is given by
$$
\Pe(M_{1}=m_{1},\ldots,M_{n}=m_{n}|\w,\X)\propto\prod_{i=1}^{n}\varphi(i,m_{i}|\w){\mbox
{ for }} \m \in \Xi.$$
\item[(iii)]The density of $\W|\MM,\X$ is
$$f(\w|\m,\X)\propto{\mbox
e}^{-\Phi_{0}(r_{1})}\prod_{i=1}^{n}{\mbox
e}^{-\Phi_{i}(w_i|r_{i+1})}\varphi(i,m_{i}|\w)\phi(X_{i}|\mu\Delta+\beta
w^{-1}_{i},w^{-1}_{i}).$$ \EndEnumerate \qed
\end{prop}
\Proof The proof is straightforward. It essentially follows from a
relabeling of the components in Theorem 4.1, combined with the
form of the likelihood in Proposition 5.2\EndProof We next
describe the posterior distribution of the integrated OU process.
\begin{prop} Suppose that $\tau$ is defined by
~\mref{model1}, then the posterior distribution of
$\tau|\MM,\W,\X$ is equivalent to the conditional distribution of
the random measure,
$$
\tau_{n}(t)=\lambda^{-1}[(1-{\mbox e}^{-\lambda t})v_{0,n}+
\int_{0}^{t}(1-{\mbox
e}^{-\lambda(t-y)})\mu_{\Omega_{n}}(dy)+\sum_{i=1}^{n}\xi_{i}{\tilde
J_{i}}(1-{\mbox e}^{-\lambda(t-{\tilde Y}_{i})})I_{\{{\tilde
Y}_{i}\leq t\}}]$$ where $v_{0,n}:=\int_{-\infty}^{0}{\mbox
e}^{y}\mu_{\Omega_{n}}(dy)+\sum_{i=1}^{n}\xi_{i}{\tilde
J_{i}}{\mbox e}^{{\tilde Y}_{i}}I_{\{{\tilde Y}_{i}\leq 0\}}$ has
the posterior distribution of $v_{0}.$\qed
\end{prop}
\Proof This result follows from Theorem 4.1 using the fact that
$\sum_{j=1}^{n(\p)}J_{j,n}\delta_{Y^{*}_{j}}\overset
{d}=\sum_{i=1}^{n}\xi_{i}{\tilde J_{i}}\delta_{{\tilde Y}_{i}}$.
See Ho~(2005) for a similar argument.\EndProof

The structure $\m\in \Xi$ is in one to one correspondence to what
are called s-paths by Brunner and Lo~(1989). See that work, in
particular Brunner and Lo~(1989, Lemma 2.1, Theorem 2.1) for
slightly different representations.  We close by noting that we
are rather surprised that the BNS-OU model used in our likelihood
structure generates s-paths. As mentioned earlier, s-paths are
known to be generated by representing monotone hazard rates as
$\int_{-\infty}^{\infty}I_{\{t<u\}}\mu(du)$ where $\mu$ is a
general completely random measure. This is the  formulation
recently investigated by~Ho~(2005) where the Gamma process results
of Dykstra and Laud~(1981), and Lo and Weng~(1989) are special
cases. Naturally, these are closely connected to the symmetric
unimodal density Dirichlet mixture models considered by Brunner
and Lo~(1989). The BNS-OU models represent the first non-trivial
{\it mixture models} outside the above mentioned class where
inference can be based on sampling solely s-paths or equivalently
$\m$, rather than $\p$. This is an important fact, since while
indeed sampling partitions $\p$ is not difficult, these models are
significantly less complex than models which can be minimally
expressed in terms of partitions. That is, the space of partitions
of the integers $\{1,\ldots,n\}$ is known to contain Bell's number
of terms which is approximately $n!$ and is considerably larger
than $\Xi$. Ho~(2005), has devised efficient computational
procedures for sampling s-path models which can be easily imported
to the present setting. See also Ho~(2002). \vskip0.2in
\centerline{\Heading References} \vskip0.2in \tenrm
\def\smc{\tensmc}
\def\sl{\tensl}
\def\bf{\tenbold}
\baselineskip0.15in \Ref  \by Barndorff-Nielsen, O. E.\yr 1997
\paper  Normal inverse Gaussian distributions and stochastic
volatility modelling\jour  Scand. J. Statist. \vol 24 \pages
1-13\EndRef

\Ref \by Barndorff-Nielsen, O.E. and Shephard, N. \yr 2001a \paper
Ornstein-Uhlenbeck-based models and some of their uses in
financial economics \jour \JRSSB \vol 63 \pages 167-241 \EndRef
\Ref \by Barndorff-Nielsen, O.E. and Shephard, N. \yr 2001b \paper
Modelling by L\'evy processes for financial econometrics. In
L\'evy processes. Theory and applications. Edited by Ole E.
Barndorff-Nielsen, Thomas Mikosch and Sidney I. Resnick. p.
283-318. Birkh\"auser Boston, Inc., Boston, MA \EndRef

\Ref \by Barndorff-Nielsen, O. E. and Shephard, N.\yr 2003 \paper
Integrated OU processes and non-Gaussian OU-based stochastic
volatility models \jour Scand. J. Statist. \vol 30 \pages 277-295
\EndRef

\Ref \by Blackwell, D. and MacQueen, J. B. \yr 1973 \paper
Ferguson distributions via P\'olya urn schemes \jour \AnnStat \vol
1 \pages 353-355 \EndRef \Ref \by Brix, A. \yr 1999 \paper
Generalized Gamma measures and shot-noise Cox processes \jour Adv.
in Appl. Probab. \vol 31 \pages 929-953 \EndRef
 \Ref \by Brunner, L.~J. and Lo, A.~Y. \yr 1989 \paper Bayes methods for a symmetric unimodal density and its
 mode \jour \AnnStat \vol 17 \pages 1550-1566 \EndRef

\Ref \by Carr, P., Geman, H., Madan, D.B. and Yor, M. \yr 2003
\paper Stochastic volatility for L\'evy processes \jour Math.
Finance \vol 13 \pages 345-382 \EndRef

\Ref \by Carr, P. and Wu, L. \yr 2004 \paper Time-changed L\'evy
processes and option pricing \jour Journal of Financial Economics
\vol 71 \pages 113-141 \EndRef

\Ref \by Daley, D. J. and Vere-Jones, D. \yr 1988 \book An
introduction to the theory of point processes \publ
Springer-Verlag \publaddr New York \EndRef

\Ref \by Duan,  J.\yr 1995 \paper The GARCH option pricing model
\jour Math. Finance \vol 5 \pages 13-32 \EndRef \Ref \by Dykstra,
R. L. and Laud, P. W. \yr 1981 \paper A Bayesian nonparametric
approach to reliability \jour \AnnStat \vol 9 \pages 356-367
\EndRef \Ref \by Eberlein, E. \yr 2001 \paper Application of
generalized hyperbolic L\'evy motions to finance. In L\'evy
processes. Theory and applications. Edited by Ole E.
Barndorff-Nielsen, Thomas Mikosch and Sidney I. Resnick. p.
319-336. Birkh\"auser Boston, Inc., Boston, MA \EndRef

\Ref \by Engle, R. F.\yr 1982 \paper Autoregressive conditional
heteroscedasticity with estimates of the variance of United
Kingdom inflation \jour Econometrica \vol 50 \pages 987-1007
\EndRef \Ref \by Escobar, M.D. \yr 1994 \paper Estimating normal
means with the Dirichlet process prior \jour \JASA \vol 89 \pages
268-277 \EndRef

\Ref \by Fiorentini, G., Sentana, E. and Shephard, N.\yr 2004
\paper Likelihood-based estimation of latent generalized ARCH
structures \jour Econometrica \vol 72 \pages 1481-1517 \EndRef

\Ref \by Ho, M.-W. \yr 2002 \paper Bayesian inference for models
with monotone densities and hazard rates. Unpublished Ph.D.
thesis. Department of Information and Systems Management. Hong
Kong University of Science and Technology\EndRef

\Ref \by Ho,M.-W.\yr 2005 \paper A Bayes method for a monotone
hazard rate via s-paths. arXiv:math.ST/0502432\EndRef

\Ref \yr 2003 \by    Ishwaran,~H. and James,~L.~F. \paper
Generalized weighted Chinese restaurant processes for species
sampling mixture models \jour  Statistica Sinica \vol   13 \pages
1211-1235 \EndRef \Ref \yr 2004 \by    Ishwaran,~H. and
James,~L.~F. \paper Computational methods for multiplicative
intensity models using
       weighted gamma processes: proportional hazards, marked point
       processes and panel count data
\jour  \JASA \vol   99 \pages 175-190 \EndRef

\Ref \by James, L.F. \yr 2002 \paper Poisson process partition
calculus with applications to exchangeable models and Bayesian
nonparametrics. arXiv:math.PR/0205093 \EndRef

\Ref \by James, L.F. \yr 2005a \paper Bayesian Poisson process
partition calculus with an application to Bayesian L\'evy moving
averages. To appear in {\it Annals of Statistics}.\\Available at
http://ihome.ust.hk/$\sim$lancelot/\EndRef

\Ref\by James, L.F. \yr 2005b \paper Analysis of a class of
likelihood based continuous time stochastic volatility models
including Ornstein-Uhlenbeck models in financial economics.
arXiv:math.ST/0503055 \EndRef

\Ref \by Jensen, M.B. and Lunde, A. \yr 2001 \paper The NIG-S and
ARCH model: a fat-tailed, stochastic, and autoregressive
conditional heteroskedastic volatility model \jour Econom. J. \vol
4 \pages 319-342\EndRef \Ref \by Korwar, R. M. and  Hollander, M.
\yr 1973 \paper Contributions to the theory of Dirichlet processes
\jour \AnnProb \pages 705-711 \EndRef

\Ref \by Liu, J.S. \yr 1996 \paper Nonparametric hierarchichal
Bayes via sequential imputations \jour \AnnStat \vol 24 \pages
911-930\EndRef \Ref \by Lo, A. Y. \yr 1984 \paper On a class of
Bayesian nonparametric estimates: I. Density Estimates \jour
\AnnStat \vol 12 \pages 351-357 \EndRef \Ref \by Lo, A. Y. and
Weng, C. S. \yr 1989 \paper On a class of Bayesian nonparametric
estimates: II. Hazard
       rates estimates
\jour Ann. Inst. Stat. Math. \vol 41 \pages 227-245 \EndRef \Ref
\by Lo, A.Y., Brunner, L.J. and Chan, A.T. \yr 1996 \paper
Weighted Chinese restaurant processes and Bayesian mixture
       model. Research Report Hong Kong University of
       Science and Technology. Available at http://www.erin.utoronto.ca/~jbrunner/papers/wcr96.pdf\EndRef
\Ref \by Prause, K.\yr 1999 \paper The generalized hyperbolic
model: estimation, financial derivatives, and risk measures. Dissertation, Institut f\"ur Mathematische Stochastik Universit\"at Freiburg \\
Available at http://www.freidok.uni-freiburg.de/volltexte/15/
\EndRef \Ref \by West, M.,  M\"uller, P. and Escobar, M.D.\yr 1994
\paper Hierarchical priors and mixture models, with application in
regression and density estimation. Aspects of uncertainty,
363-386, Wiley Ser. Probab. Math. Statist. Probab. Math. Statist.,
Wiley, Chichester, 1994 \EndRef \smc

\Tabular{ll}

Lancelot F. James\\
The Hong Kong University of Science and Technology\\
Department of Information and Systems Management\\
Clear Water Bay, Kowloon\\
Hong Kong\\
\rm lancelot\at ust.hk\\

\EndTabular
\medskip

\Tabular{ll}

John W. Lau\\
Department of Mathematics\\
University of Bristol\\
Bristol, BS8 1TW, UK\\
\rm John.Lau\at bristol.ac.uk \\
\EndTabular
\end{document}